\documentclass{amsart}
\usepackage{pstricks}
\usepackage{pst-all}
\usepackage{amscd}
\usepackage{stmaryrd}
\usepackage{amssymb}
\usepackage{enumerate}
\usepackage{amsmath}
\usepackage{amsthm}
\usepackage{graphicx}
\usepackage[all,cmtip]{xy}
\usepackage{geometry}

\DeclareMathAlphabet{\mathpzc}{OT1}{pzc}{m}{it}

\theoremstyle{definition}

\DeclareMathOperator{\inv}{inv}
\DeclareMathOperator{\ord}{ord}

\newcommand{\N}{\mathbb{N}}

\newcommand{\po}{\frac{p-1}{2}}
\newcommand{\ut}{A}
\newcommand{\vt}{B}
\newcommand{\ot}{Q}
\newcommand{\y}{u}
\newcommand{\K}{\mathbb {K}}
\newcommand{\T}{\Delta}
\newcommand{\TS}{\Gamma}
\newcommand{\TT}{B}
\newcommand{\clean}{_{\textnormal{clean}}}
\newcommand{\st}{} 

\font\tss=ptmr at 8pt

\long\def\ignore#1\recognize{}

\def\weak{strict }
\def\A{{\mathbb A}}

\begin{document}

\title[Cycles of singularities]{Cycles of singularities appearing in the resolution problem in positive characteristic}
\author{Herwig Hauser, Stefan Perlega}
\maketitle

\let\thefootnote\relax\footnote{
{\tss \hskip -.41cm
MSC-2000: 14B05, 14E15, 12D10.  We are grateful to an anonymous referee for valuable suggestions 
improving the presentation of the article. Supported by project P-25652 of the Austrian Science Fund FWF. }}

\begin{abstract} We present a hypersurface singularity in positive characteristic which is defined by a purely inseparable power series, and a sequence of point blowups so that, after applying the blowups to the singularity, the same type of singularity reappears after the last blowup, with just certain exponents of the defining power series shifted upwards. The construction hence yields a cycle. Iterating this cycle leads to an infinite increase of the residual order of the defining power series. This disproves a theorem claimed by Moh about the stability of the residual order under sequences of blowups. It is not a counter-example to the resolution in positive characteristic since larger centers are also permissible and prevent the phenomenon from happening.



\end{abstract}



\def\mult{{\rm mult}}
\def\qq{q}
\def\aa{a}
\def\bb{b}
\def\rr{r}
\def\ss{s}
\def\kk{k}


\section{Introduction}

\noindent The Whitney-umbrella is the surface $X$ in $\A^3$ defined by the equation 
\[x^2-y^2z=0.\]
 Its singular locus is the $z$-axis, and $0$ is the worst singularity; at the other points of the $z$-axis, the singularity has normal crossings and consists of two smooth transversal branches.  A simple computation shows that blowing up $\A^3$ with center the origin $0$ produces a strict transform $X'$ of $X$ which has at the origin of the $z$-chart the same singularity as $X$ at $0$. Thus a cycle occurs, and no improvement of the singularity is achieved. It turns out that the appropriate choice of center is the whole $z$-axis in $\A^3$: this blowup resolves the singularities of $X$ in one stroke.

Consider more generally a surface in $\A^3$ of equation $x^a-y^bz^c=0$ with $a<b,c$. Blowing up the origin produces in the $z$-chart the equation $x^a-y^bz^{c+b-a}=0$. As $c'=c+b-a>c$ the singularity has gotten worse, while the shape of the equation is the same. This type of repetition will be called a {\it cycle under blowup}, with {\it shifted} exponents (the $z$-exponent $c$ transforms into $c'$, the others remain constant).

In characteristic zero, the occurrence of cycles under blowup can always be avoided by the correct choice of the center: There is a systematic way via hypersurfaces of maximal contact and coefficient ideals to determine an appropriate (regular) center so that the induced blowup does improve the singularity \cite{Hironaka_Annals, Villamayor_Constructiveness, Villamayor_Patching, BM_Canonical_Desing, Wlodarczyk, EH}. This is the way to prove resolution of singularities in zero characteristic (up to a bunch of technicalities).

The situation is much more subtle in the case of positive characteristic $p>0$. We will produce, for any $p$, examples where cycles of singularities appear after a sequence of blowups in a much more involved fashion. The examples exhibit a new difficulty which adds to the many already known problems one encounters in positive characteristic. The blowups in the examples always have point centers. 
Taking larger centers would prevent the phenomenon from happening. We were not able to construct examples with cycles where the choice of point centers is forced (e.g., because the singularities are isolated). Of course, such an example would disprove the existence of resolution of singularities in positive characteristic.

The examples are of a relatively simple form, namely purely inseparable power series of the form
\[f=z^{p^e}+x^ry^sw^t(w^dy^b+x^av^c) + \ldots,\]
up to higher order terms. The exponents and the sequence of blowups have to be chosen very carefully to produce a cycle. After all these blowups, the transformed equation has the same shape, with certain exponents shifted {\it upwards}. 
\medskip

The rest of this section recalls the main ingredients of resolution proofs which are used in the discussion of the examples; it can be skipped by the expert reader. 

Let $X$ be a singular variety embedded in a regular ambient variety $W$. For every point $a$ of $X$, one may choose locally at $a$ a defining ideal $J$ for $X$. The order $\ord_aX:=\ord_aJ$ of $J$ at $a$ is one of the main numerical ingredients for the construction of the resolution of the singularities of $X$: it stratifies $X$ into a disjoint union of finitely many locally closed strata, the smallest one collecting the worst singularities of $X$ (we assume here that the embedding of $X$ in $W$ is minimal). Blowing up regular centers inside this stratum, one then aims at lowering the maximum of the local orders by a suitable sequence of blowups. The key observation here is that this maximum cannot increase under such a {\it permissible} blowup. If it drops, one may apply induction to achieve a resolution. If it remains constant, an extra argument is necessary to proceed. Typically, one then tries to apply induction on the ambient dimension: Choosing suitable regular hypersurfaces $V$ in $W$ locally at the points of the worst stratum, one associates to the variety $X$ ideals $K$ in $V$ whose complexity reflects the original singularity and whose resolution simplifies sufficiently $X$ so as to allow likewise its resolution. 

In characteristic zero, this program works due to the existence of hypersurfaces of {\it maximal contact}. They ensure that all local constructions patch and that the local descent in dimension to the ideals $K$ commutes with blowup at those points of the exceptional divisor where the multiplicity of $X$ has remained constant.

For the inductive argument it is necessary to factor the local ideals $K$ into products $K=M\cdot I$, where $M$ is a principal monomial ideal defining the exceptional divisor produced by the earlier blowups in the resolution process. The residual factor $I$ carries the relevant information on the complexity of $X$ locally at the respective point of the stratum. Its order $\ord_aI$ forms the second component of a local upper semicontinuous invariant $\inv_aX$ of $X$ at $a$. This invariant consists of a string of orders in descending dimensions,
\[\inv_aX=(\ord_aX,\ord_aI, ...),\]
which then, taken lexicographically, serves as a resolution invariant for $X$ (leaving aside technicalities): Namely, it refines the stratification of $X$ by the order of $J$, yielding a smallest stratum which is regular, and, so that, when blowing up  this stratum, the invariant drops at every point of the exceptional divisor. Then induction applies to establish resolution.

In positive characteristic, this approach meets serious obstructions: First, hypersurfaces of maximal contact need no longer exist. A substitute are regular hypersurfaces maximizing the order of the ideals $K$ and $I$ (which is also the case for hypersurfaces of maximal contact). They permit the definition of $\ord_aI$ as an intrinsic and significant second component of $\inv_aX$, the {\it residual order} of $X$ at $a$. The next obstruction then is the fact that in positive characteristic the residual order may increase under blowup at points where $\ord_aX$ has remained constant. However, the increase is not too large, happens very rarely, and can explictly be bounded, as Moh showed \cite{Moh, Ha_BAMS_2, HP_Characterizing}. This raised the hope that the invariant $\inv_aX$ may drop {\it in the long run} and thus serves again to establish resolution. This has been proven to work for surfaces \cite{HW, HP_Surfaces} and is still open in higher dimensions. For arbitrary dimensions, Moh claimed in \cite{Moh} that the invariant cannot increase under permissible blowup beyond a certain bound. It turns out that this claim is false. In fact, we prove in the present note:\medskip

{\it There exists a purely inseparable hypersurface singularity $X$ and a sequence of permissible blowups along which the order of $X$ remains constant but for which the orders of the residual ideals $I$ tend to infinity.}\medskip

The hypersurface and the sequence of blowups were constructed by the second author (preceded by an earlier and less significant example of the first author). The centers are always points, chosen carefully inside the respective exceptional divisor. After a tricky composition of such blowups the defining equation of $X$ will have transformed into one of exactly the same shape, but with certain exponents and the residual order $\ord_aI$ shifted upwards. The iteration of the blowups then makes the residual order go to infinity.




\section{Setting}


\noindent Let $\K$ be an algebraically closed field of characteristic $p>0$. The kind of singularities that we consider in this article are defined by \emph{purely inseparable equations} of the form
\[f(z,x)=z^{p^e}+F(x_1,\ldots,x_n)=0,\]
where $e$ is a positive integer and $F\in \K[[x_1,\ldots,x_n]]$ is a power series \st of order $\ord F\geq p^e$. 


The key observation here 
is that a change of parameters $z_1=z-g$ \st with $g\in \K[[x_1,\ldots,x_n]]$ changes the expansion of $f$ to
\[f=z_1^{p^e}+g(x_1,\ldots,x_n)^{p^e}+F(x_1,\ldots,x_n).\]
Hence, any $p^e$-th power that appears in the expansion of $F$ can be eliminated via such a change.
This suggests to consider all $p^e$-th powers in the power series expansion of $F$ as artificial information. We then say that $F$ is \emph{clean} if no $p^e$-th powers appear in its expansion. Respectively, we refer to a change of parameters $z_1=z-g$ which eliminates all $p^e$-th powers from 
$F$ as \emph{cleaning}.

    
To measure the improvement of the singularity defined by $f$ under blowup, it is standard to make use of a normal crossings divisor $E$ of the form $E=V(\prod_{i\in\T}x_i)$ for some subset $\T\subseteq\{1,\ldots,n\}$. This divisor consists of the exceptional components produced by earlier blowups in the resolution process.

    
Consider now a purely inseparable power series $f=z^{p^e}+F(x_1,\ldots,x_n)$ where $F$ is clean, together with  $E=V(\prod_{i\in\T}x_i)$. Write $F$ in the form
\[F(x_1,\ldots,x_n)=\prod_{i\in\Delta}x_i^{r_i}\cdot G(x_1,\ldots,x_n)\]
where $G\in\K[[x_1,\ldots,x_n]]$ is a power series and $r_i=\ord_{(x_i)}F$ denotes the order of $F$ along the component $x_i=0$. We call the order of this power series $G$ the \emph{residual order} of $f$ with respect to the normal crossings divisor $E$,
\[{\rm residual.order }_E\, f =\ord\, G.\]
 This is a natural invariant which measures how far, up to multiplication by units in the power series ring, $F$ is away from being a monomial in $x_1,\ldots,x_n$. The residual order and similar numerical characters as well as the cleaning process have appeared repeatedly, often only implicitly, in the literature on resolution of singularities in positive characteristic \cite {Abhyankar_67, Moh, Cossart_Piltant_1, Cossart_Piltant_2, Cutkosky_Skeleton, BV_Monoidal, Kawanoue_Matsuki_Surfaces}. The terminology residual order was proposed by Hironaka in \cite{Hironaka_CMI}. 
%



 Notice further that, if $r_i\geq p^e$, we may blow up the center $z=x_i=0$ and consider the transform of $f$ in the $x_i$-chart. This reduces $r_i$ by $p^e$ and leaves the power series $F$ otherwise unchanged.

 One can view the residual order as a generalization of a resolution invariant that is successfully used over fields of characteristic zero. 
We refer to \cite{EH, Ha_BAMS_2, HP_Characterizing} for this interpretation.


An ideal $P$ of the ring $\K[[z,x_1,\ldots,x_n]]$ is said to define a \emph{permissible center} of blowup for $f$ and $E$ as above if $P$ is of the form $P=(z,x_i:i\in\TS)$ for some subset $\Gamma\subseteq\{1,\ldots,n\}$ and the following two conditions hold:
\begin{enumerate}

 \item $f\in P^{p^e}$.
 \item $G\in P^d$, where $G$ is defined via the factorization $F=\prod_{i\in\Delta}x_i^{r_i}\cdot G$ as above and $d$ is the residual order of $f$.

\end{enumerate}

\noindent The maximal ideal $(z,x_1,\ldots,x_n)$ always defines a permissible center.

    
Now consider a blowup map $\pi:\K[[z,x_1,\ldots,x_n]]\to\K[[z,x_1,\ldots,x_n]]$ with center $P$ as above and let $f'$ denote the \weak transform of $f$. It is well-known that the inequality
\[\ord f'\leq\ord f=p^e\]
always holds. If the order of $f'$ is strictly smaller than $p^e$, the singularity has already significantly improved under this blowup \st 
. Hence, we only consider the case where equality holds. There then exists an index $j\in\Gamma$ and constants $t_i\in\K$ such that the map $\pi$ is of the form
\[\begin{array}{ll}
 \pi(z)=x_jz,&\\
 \pi(x_j)=x_j, &\\ 
 \pi(x_i)=x_j(x_i+t_i) & \text{for $i\in \TS\setminus\{j\}$,}\\
 \pi(x_l)=x_l & \text{for $l\notin \TS$.}
\end{array}\]
If $P$ is the maximal ideal, we will refer to $\pi$ as a \emph{point blowup}.

For the examples below, we adapt the following terminology for point blowups: Only the parameter $x_j$ (referred to as the $x_j$-chart) and the translations $x_i\mapsto x_i+t_i$ with non-zero constants $t_i$ are mentioned. For instance, the expression
 \[\text{$x_1$-chart, $x_2\mapsto x_2+1$}\]
 encodes the map $\pi$ given by $\pi(z)=x_1z$, $\pi(x_1)=x_1$, $\pi(x_2)=x_1(x_2+1)$ and $\pi(x_i)=x_1x_i$ for $i\geq 3$.


The \weak transform of $f$ under $\pi$ has the form
\[f'=z^{p^e}+F'(x_1,\ldots,x_n)\]
where $F'=x_j^{-p^e}\pi(F)$. Notice that, even though $F$ is clean, the power series $F'$ is not necessarily clean again. After applying cleaning, we obtain a new expansion
\[f'=z^{p^e}+F'\clean(x_1,\ldots,x_n)\]
from which we can compute the residual order of $f'$ with respect to the normal crossings divisor $E'=V(x_j\cdot\prod_{i\in\Delta\setminus \TT}x_i)$ \st obtained as the total transform of $E$, where the set $\TT$ is defined as $\TT=\{j\}\cup\{i:t_i\neq0\}$.

\section{Behaviour of the residual order under blowup}

\noindent Contrary to the behavior of the resolution invariants used in characteristic zero, it is known that the residual order of a purely inseparable power series $f=z^{p^e}+F(x_1,\ldots,x_n)$ may increase under permissible blowups. Moh was able to show in \cite{Moh} that the increase under a \emph{single blowup} is bounded by $p^{e-1}$. A further analysis of the increase and necessary conditions that $F$ has to fulfill for an increase to happen were given by the first author in \cite{Ha_BAMS_2}, see also \cite{HP_Characterizing}.

Moh claimed in \cite{Moh} that the following \emph{Stability Theorem} holds: \medskip

\noindent \emph{Let be given a purely inseparable equation $f=z^{p^e}+F(x_1,\ldots,x_n)=0$ with residual order $d$ and a sequence of permissible blowups under which the order of the \weak transforms of $f$ remain constant. Then the residual orders of the \weak transforms of $f$ cannot increase beyond the bound $d+p^{e-1}$.} \medskip

In particular, this result would rule out the possibility of the residual order increasing indefinitely.

Below we will present examples of purely inseparable equations  with $e=3$ and infinite sequences of point blowups under which the residual order tends to infinity. This disproves Moh's claim in the case $e\geq 3$ (it is known to be valid for $e=1$).

Although the examples might seem discouraging for proving resolution of singularities in positive characteristic, we emphasize that they do not constitute a counterexample to the existence of resolutions in positive characteristic. The reason for this is that in the examples one could choose at various instances a larger center than a point and thus would end up with a different sequence of blowups for which the residual order need not tend to infinity.

The flaw in Moh's proof of the Stability Theorem is relatively subtle. On p.~970 in \cite{Moh}, a non-negative integer $r$ is defined as the maximal number with the property that the initial form of $F$ is a $p^r$-th power. It is then claimed that one may assume without loss of generality that the residual order $d$ is divisible by $p^r$. This leads to the false conclusion that, after blowup, the residual order is bounded by $d'\leq[d/p^r]\cdot p^r+p^r$. As the following example shows, this is not true. Consider the purely inseparable polynomial
\[f=z^4+x^2y^2w^3(w(x+y)^4+x^{13})\]

\noindent over a field of characteristic $2$. The residual order of $f$ is $d=5$. Further, $p^r=2$. Let $f'$ be the \weak transform of $f$ under a point blowup in the $x$-chart with translation $y\mapsto y+1$. Then
\[\begin{array}{l}
   f'=z^4+x^8w^3(wy^4+wy^6+x^8+x^8y^2)
    =z_1^4+x^8w^3(wy^6+x^8+x^8y^2)
  \end{array}
\]

\noindent after cleaning. Hence, the residual order of $f'$ is $d'=7$. This exceeds the bound $[d/p^r]\cdot p^r+p^r=6$.


\section{Examples for cycles and the indefinite increase of the residual order}

We provide two types of examples, the first covering characteristic $2$ and the second all odd characteristics:  After a well chosen sequence of point blowups, the defining equation has the same shape as at the beginning, but some exponents have increased. Also, the residual order has increased.

In the examples, the letters $\ut$ and $\vt$, respectively $\ot$, will denote unspecified units $\ut,\vt \in \K[[x_1,\ldots,x_n]]^*$, respectively power series $\ot\in \K[[x_1,\ldots,x_n]]$, and $\lambda$ will denote an unspecified non-zero constant $\lambda\in \K^*$. These objects can be chosen arbitrarily at the beginning of the sequence. Under the blowups, the objects maintain their quality but change their value according to the transformation rules dictated by the blowups. To keep the formulas readable, the same letters are used after each blowup and when rewriting the polynomials. 

The consecutive \weak transforms of the purely inseparable powers series $f$ under the sequence of blowups will all be denoted again by $f$.\goodbreak

\medskip


\noindent {\it First example}: Field of characteristic $p=2$, $n=5$, $\ord f = 8$.\medskip 

\noindent The variables $x_1,\ldots,x_5$ will be denoted by $x,y,u,v,w$.
Choose for $d$ any fixed {\it even} positive integer. Let $\aa$, $\bb$, $\rr$ and $\ss$ be non-negative integers chosen arbitrarily at the beginning, with $8\ss \geq d$. These latter numbers specify the exceptional multiplicities at the start. In the course of the blowups, the values of $\aa$, $\bb$, $\rr$ and $\ss$ change, though will  be denoted by the same letters to avoid complicated formulas. The transformation rules under blowup ensure that the exponents never become negative. As $\aa$, $\bb$, $\rr$ and $\ss$ are always multiplied by $8=\ord\, f$, their actual values do not affect the validity of the cleaning process. 
One could, alternatively, apply auxiliary blowups with codimension $2$ centers of the type $(z,x)$, $(z,y)$, etc., considered always in the $x$-chart, $y$-chart, etc., in order to reduce after each step the exceptional exponents modulo $8$.

All blowups are point blowups. 

\begin{enumerate}  \setcounter{enumi}{-1}

 
 \item Starting equation:
 

$f=z^8+x^{8\aa+4}y^{8\bb+4}u^{8\rr}w^{8\ss-d}\cdot(w^d(\lambda+u^{2d+6}\cdot\ot)+x^{d+1}\y^{2d+6}\cdot\ut).$

The residual order equals $d$.


 \item $x$-chart, monomial blowup:
 
 $f=z^8+x^{8\aa}y^{8\bb+4}u^{8\rr}w^{8\ss-d}\cdot(w^d(\lambda
+x^{2d+6}\cdot\ot)+x^{2d+7}\cdot\ut)$

 Translation $y\mapsto y+1$, $u\mapsto u+1$:
 
$f=z^8+x^{8\aa}(y+1)^{8\bb+4}(u+1)^{8\rr}w^{8\ss-d}\cdot(w^d(\lambda
+x^{2d+6}\cdot\ot)+x^{2d+7}\cdot\ut)$

\hskip .3cm  $=z^8+x^{8\aa}w^{8\ss-d}\cdot(\lambda (u+1)^{8\rr}w^d+w^dy^4\cdot\vt+w^dx^{2d+6}\cdot\ot+x^{2d+7}\cdot\ut),$

after multiplying $(y+1)^{8\bb+4}(u+1)^{8\rr}$ with the residual terms in the parenthesis, getting 

new invertible series $A,B$ and a new arbitrary series $Q$ on the second line. Cleaning gives

$f=z_1^8+x^{8\aa}w^{8\ss-d}\cdot(w^dy^4\cdot\vt+w^dx^{2d+6}\cdot\ot+x^{2d+7}\cdot\ut).$

\st At this point, the residual order has increased to $d+4$.


 \item $x$-chart: \st
 
$f=z^8+x^{8\aa+4}w^{8\ss-d}\cdot(w^dy^4\cdot\vt+w^dx^{2d+2}\cdot\ot+x^{d+3}\cdot\ut)$

\hskip .3cm $=z^8+x^{8\aa+4}w^{8\ss-d}\cdot(w^dy^4\cdot\vt+x^{d+3}\cdot\ut).$

The residual order has decreased to $d+3$. 


 \item $\y$-chart:
 
$f=z^8+x^{8\aa+4}u^{8\rr+7}w^{8\ss-d}\cdot(w^dy^4u\cdot\vt+x^{d+3}\cdot\ut).$


\item $\frac{d}{2}$ times $v$-chart:

$f=z^8+x^{8\aa+4}u^{8\rr+7}v^{(8\aa+8\rr+8\ss+6)d/2}w^{8\ss-d}\cdot(w^dy^4uv^{d}\cdot\vt+x^{d+3}\cdot\ut).$


\item $w$-chart, $y\mapsto y+1$, $v\mapsto v+1$:

$f=z^8+x^{8\aa+4}u^{8\rr+7}(v+1)^{(8\aa+8\rr+8\ss+6)d/2}w^{8\ss-(d+2)}\cdot(w^{d+2}u\cdot\vt+x^{d+3}\cdot\ut)$

\hskip .3cm $=z^8+x^{8\aa+4}u^{8\rr+7}w^{8\ss-(d+2)}\cdot(w^{d+2}u\cdot\vt+x^{d+3}\cdot\ut)$,

\st after multiplying $(v+1)^{(8\aa+8\rr+8\ss+6)d/2}$ with the residual terms in the parenthesis. 


\item $u$-chart:

$f=z^8+x^{8\aa+4}u^{8\rr+4}w^{8\ss-(d+2)}\cdot(w^{d+2}\cdot\vt+x^{d+3}\cdot\ut).$

The residual order has decreased to $d+2$.


\item $(2d+10)$ times $y$-chart:

$f=z^8+x^{8\aa+4}y^{8\bb}u^{8\rr+4}w^{8\ss-(d+2)}\cdot(w^{d+2}(\lambda+y^{2d+10}\cdot\ot)+x^{d+3}y^{2d+10}\cdot\ut).$

The residual order equals $d+2$.

\end{enumerate}


\noindent At this point, $f$ has again the same form as the starting equation (exchanging $y$ with $u$ and $\bb$ with $\rr$), with the exception of $d$ being raised to $d+2$. The residual order has increased by $2$ during the above sequence of blowups. The final value of $\ss$ is such that the exponent $8\ss-(d+2)$ of $w$ in the last equation is again non-negative. Since the sequence of blowups can be iterated indefinitely, the residual order also increases indefinitely.
\medskip\goodbreak

The sequence of blowups has a relatively simple structure: The residual order increases only under blowup (1). The remaining blowups in the sequence are necessary to rebuild the starting equation with increased $d$. After blowup (2), the residual order decreases by $1$ and the initial form of $G$ changes to $x^{d+3}$. The blowup sequences (3) and (4) leave the initial form unchanged, but add powers of $u$ and $v$ to the term $w^{d}y^4$. These powers are necessary to create the new term $w^{d+2}u\cdot\vt$ with blowup (5). Under blowup (6) the residual order decreases again by $1$ and $w^{d+2}$ becomes the new initial form of $G$. At this point, the initial form already has the same form as in the starting equation. The last sequence of blowups adds powers of $y$ to $x^{d+3}$ until it is of the same form as in the starting equation. It also has the effect of transforming the unit $V$ into a constant $\lambda$ plus a higher order term that is divisible by $y^{2d+10}$.



\bigskip\goodbreak


\noindent {\it Second example}: Field of characteristic $p\geq 3$, $n=4$, $\ord f= p^3$. \medskip

\noindent The parameters $x_1,\ldots,x_4$ will be denoted by $x,y,v,w$. Choose for $d$ a fixed positive integer that is divisible by $p\po $. 
Let $\aa$, $\bb$, $\rr$ and $\ss$ be integers chosen arbitrarily at the beginning but subject to the inequalities $\aa\geq 0$, $\bb \cdot p^3\geq \frac{1}{2}(p^3-1)$, $\rr\geq 0$ and $\ss \cdot p^3\geq d$. These numbers specify the exceptional multiplicities at the start. The inequalities ensure that all exponents of the exceptional monomial in the starting equation are non-negative. In the course of the blowups, the values of $\aa$, $\bb$, $\rr$ and $\ss$ change, though will  be denoted by the same letters. The transformation rules under blowup imply that the exponents never become negative. Moreover, as $\aa$, $\bb$, $\rr$ and $\ss$ are always multiplied by the order $p^3$ of $f$, their actual values do not affect the validity of the cleaning process. Again, one could alternatively apply auxiliary blowups with codimension $2$ centers in order to reduce after each step the exceptional exponents modulo $p^3$.

Set $d'=d+\po p$, $m=\frac{2d}{p-1}+p-1$, and $\qq={\frac{p+1}{2}(d+p^2-1)}$, $\qq'={\frac{p+1}{2}(d'+p^2-1)}$. These are fixed values which do not change in the formulas below. All blowups are point blowups, $A$ and $B$ denote again unspecified invertible series, $Q$ an arbitrary series, and $\lambda\in \K^*$ a constant. 


\begin{enumerate}  \setcounter{enumi}{-1}
 
 \item Starting equation:


$f=z^{p^3}+x^{\aa p^3+\frac{p^3-p^2}{2}}y^{\bb p^3-\frac{p^3-1}{2}}v^{\rr p^3}w^{\ss p^3-d}\cdot(y^{\frac{p^2-1}{2}}w^d(\lambda+v^\qq\cdot\ot)+ x^{d+\frac{p^2+1}{2}}v^\qq\cdot\ut).$

The residual order equals $\frac{p^2-1}{2}+d$. 


   
 \item $x$-chart, monomial blowup (with varying $a$, $Q$ and $A$):
 
 $f=z^{p^3}+x^{\aa p^3-\frac{p^2-1}{2}-d}y^{bp^3- \frac{p^3-1}{2}}v^{rp^3}w^{\ss p^3-d} (x^{\frac{p^2-1}{2}+d} (y^{\frac{p^2-1}{2}}w^d(\lambda+x^qv^\qq\cdot\ot)+ x^{q+1}v^\qq\cdot\ut))$\medskip

\hskip .3cm $=z^{p^3}+x^{\aa p^3}w^{\ss p^3-d}  (y^{bp^3 + \frac{1-p}{2}p^2}v^{rp^3}w^d(\lambda+x^qv^\qq\cdot\ot)+ x^{q+1}y^{bp^3- \frac{p^3-1}{2}}v^{rp^3+\qq}\cdot\ut)$\medskip

\hskip .3cm $=z^{p^3}+x^{\aa p^3}w^{\ss p^3-d}  (y^{bp^3 + \frac{1-p}{2}p^2}v^{rp^3}w^d(\lambda+x^q\cdot\ot)+ x^{q+1}y^{bp^3- \frac{p^3-1}{2}}v^{rp^3+\qq}\cdot\ut)$\medskip

Translation $y\mapsto y+1$, $v\mapsto v+1$ (with varying $Q$, $A$ and $B$):\medskip
 
 $f=z^{p^3}+x^{\aa p^3}w^{\ss p^3-d}  ((y+1)^{bp^3 + \frac{1-p}{2}p^2}(v+1)^{rp^3}w^d(\lambda+x^q\cdot\ot)+ x^{q+1}\cdot\ut)$\medskip
 
 
 \hskip .3cm $=z^{p^3}+x^{\aa p^3}w^{\ss p^3-d}(\lambda w^d(v+1)^{rp^3}+y^{p^2}w^d\cdot \vt+x^qw^d\cdot\ot+x^{q+1}\cdot\ut).$\medskip
 
Cleaning gives \medskip

$f=z_1^{p^3}+x^{\aa p^3}w^{\ss p^3-d}(y^{p^2}w^d\cdot\vt+x^\qq w^d\cdot\ot+x^{\qq+1}\cdot\ut).$\medskip
 

\ignore

From version 3: (parenthesis in first line is not correct, and second line misses $v^{rp^3} w^d\cdot\vt'$)\medskip

$f=z^{p^3}+x^{\aa p^3}(y+1)^{\bb p^3-\frac{p^3-1}{2}}(v+1)^{\rr p^3}w^{\ss p^3-d}\cdot (\lambda w^d+w^dy^{p^2}\cdot\vt+w^dx^\qq\cdot\ot+x^{\qq+1}\cdot\ut)$

\hskip .3cm $=z^{p^3}+x^{\aa p^3}w^{\ss p^3-d}\cdot (\lambda w^d+w^dy^{p^2}\cdot\vt+w^dx^\qq\cdot\ot+x^{\qq+1}\cdot\ut),$

after multiplying $(y+1)^{\bb p^3-\frac{p^3-1}{2}}(v+1)^{\rr p^3}$ with the residual terms in the parenthesis, 

getting new invertible series $A,B$ and a new arbitrary series $Q$ on the second line. 

Cleaning gives 

$f=z_1^{p^3}+x^{\aa p^3}w^{\ss p^3-d}(w^dy^{p^2}\cdot\vt+w^dx^\qq\cdot\ot+x^{\qq+1}\cdot\ut).$

\recognize

At this point, the residual order has increased by $\frac{p^2+1}{2}$ from $\frac{p^2-1}{2}+d$ to $p^2+d$.

 
\item $\po$ times $x$-chart: \st
 
$f=z^{p^3}+x^{\aa p^3+\frac{p^3-p^2}{2}}w^{sp^3-d}\cdot (y^{p^2}w^d\cdot\vt+x^{\frac{p+1}{2}(d-1)+p^2}w^d\cdot\ot+x^{d+p^2-\po}\cdot\ut)$

\hskip .3cm $=z^{p^3}+x^{\aa p^3+\frac{p^3-p^2}{2}}w^{sp^3-d}(y^{p^2}w^d\cdot\vt+x^{d'+\frac{p^2+1}{2}}\cdot\ut)$.

The residual order has decreased to $d'+\frac{p^2+1}{2}=p^2-\po+d$.


 \item $m$ times $v$-chart:
 
$f=z^{p^3}+x^{\aa p^3+\frac{p^3-p^2}{2}}v^{\rr p^3+m\frac{p^3+p^2-p+1}{2}}w^{sp^3-d}\cdot(y^{p^2}v^{m\po}w^d\cdot\vt+x^{d'+\frac{p^2+1}{2}}\cdot \ut).$


\item $w$-chart, $v\mapsto v+1$ (with varying $s$, $A$ and $B$): 

$f=z^{p^3}+x^{\aa p^3+\frac{p^3-p^2}{2}}(v+1)^{\rr p^3+m\frac{p^3+p^2-p+1}{2}}w^{sp^3-d'}\cdot (y^{p^2}w^{d'}\cdot\vt+x^{d'+\frac{p^2+1}{2}}\cdot\ut)$

\hskip.3cm $=z^{p^3}+x^{\aa p^3+\frac{p^3-p^2}{2}}w^{sp^3-d'}\cdot (y^{p^2}w^{d'}\cdot\vt+x^{d'+\frac{p^2+1}{2}}\cdot\ut)$,

after multiplying $(v+1)^{\rr p^3+m\frac{p^3+p^2-p+1}{2}}$ with the residual terms in the parenthesis. 



\item $y$-chart:

$f=z^{p^3}+x^{\aa p^3+\frac{p^3-p^2}{2}}y^{\bb p^3-\frac{p^3-1}{2}}w^{\ss p^3-d'}\cdot(y^{\frac{p^2-1}{2}}w^{d'}\cdot\vt+x^{d'+\frac{p^2+1}{2}}\cdot\ut)$.


The residual order has decreased to $\frac{p^2-1}{2}+d'=p^2-\frac{p+1}{2}+d$.


\item $\qq'$ times $v$-chart:

$f=z^{p^3}+x^{\aa p^3+\frac{p^3-p^2}{2}}y^{\bb p^3-\frac{p^3-1}{2}}v^{\rr p^3}w^{\ss p^3-d'}\cdot(y^{\frac{p^2-1}{2}}w^{d'}(\lambda+v^{\qq'}\cdot\ot) +x^{d'+\frac{p^2+1}{2}}v^{\qq'}\cdot\ut).$

The residual order equals $\frac{p^2-1}{2}+d'$. 

 

\end{enumerate}


\noindent  At this point, $f$ has again the same form as the starting equation, with the exception of $d$ being raised to $d'=d+\po p$, $q$ being replaced by $q'$, and new values for $\aa$, $\bb$, $\rr$ and $\ss$.







The residual order has increased along the sequence by $\po p$ from $d+\frac{p^2-1}{2}$ to $d'+\frac{p^2-1}{2}$. The final value of $\bb$ and $\ss$ is so that the exponents $\bb\cdot p^3-\frac{p^3-1}{2}$ of $y$ and $\ss\cdot p^3-d'$ of $w$ in the last equation are again non-negative.  Since the process can be iterated indefinitely, the residual order also increases indefinitely.


\ignore

\medskip\hrule \medskip

\noindent (*) Computation of the $w$-exponent in blowup (4): After blowup (3), we had\medskip

$f=z^{p^3}+x^{\aa p^3+\frac{p^3-p^2}{2}}v^{\rr p^3+m\frac{p^3+p^2-p+1}{2}}w^{sp^3-d}\cdot(y^{p^2}v^{m\po}w^d\cdot\vt+x^{d'+\frac{p^2+1}{2}}\cdot \ut),$\medskip

\noindent where the two terms in the parenthesis have order $p^2+m\po+d$ and $d'+\frac{p^2+1}{2}$, respectively, using the fact that $\vt$ and $\ut$ are invertible series. As $m=\frac{2d}{p-1}+p-1$ and $d'=d+p\po$ the second order is smaller than the first one and thus appears as a summand in the computation of the $w$-exponent under the next blowup. The value of the $w$-exponent in the exceptional factor of the strict transform $f$ in the first displayed formula of (4) is therefore\medskip

$w$-exp${}=ap^3+\frac{p^3-p^2}{2}+rp^3+m\frac{p^3+p^2-p+1}{2} +sp^3-d+d'+\frac{p^2+1}{2}-p^3$.\medskip

\noindent By assumption, $d$ is divisible by $p\po$, so let us write $d=c\cdot p\po$ with $c\in\N$, $c\geq 1$. This gives $m= (c+1)p-1$ and $d'=(c+1)p\po$. Substitution yields
\medskip

$w$-exp${}=(a+r+s-1)p^3+\frac{p^3-p^2}{2}+((c+1)p-1)\frac{p^3+p^2-p+1}{2} +p\po+\frac{p^2+1}{2}$,\medskip

\hskip 1cm $=(a+r+s-1)p^3+\frac{p^3+p^2-p+1}{2} +((c+1)p-1)\frac{p^3+p^2-p+1}{2} $,\medskip

\hskip 1cm $=(a+r+s-1)p^3+(c+1)p^3\frac{p+1}{2} -(c+1)p\po$,\medskip

\hskip 1cm $=(a+r+s-1+(c+1)\frac{p+1}{2})p^3-d'=:\widetilde sp^3-d'$.\medskip

\noindent This is the claimed value of the $w$-exponent, the new $s$ being $\widetilde s=a+r+s-1+(c+1)\frac{p+1}{2}$. By construction, $\widetilde sp^3-d'$ is non-negative, and we write it in (4) again as $sp^3-d'$. Notice here that the new series $\ut$ is again invertible.\medskip

\noindent Let us also check the value $d'$ of the $w$-exponent of the first term $y^{p^2}w^{d'}\cdot\vt$ inside the parenthesis appearing in the expansion of $f$ in (4).  It is given by \medskip

$p^2+m\po +d-(d'+\frac{p^2+1}{2})=p^2+((c+1)p-1)\po-p\po-\frac{p^2+1}{2}=(c+1)p\po$,\medskip

\noindent which equals $d'$ as claimed. Moreover, observe that the new series $\vt$ is again invertible. This establishes the formulas in (4).
\medskip \hrule\medskip

\recognize

\hskip .5cm \includegraphics[width=0.7\textwidth]{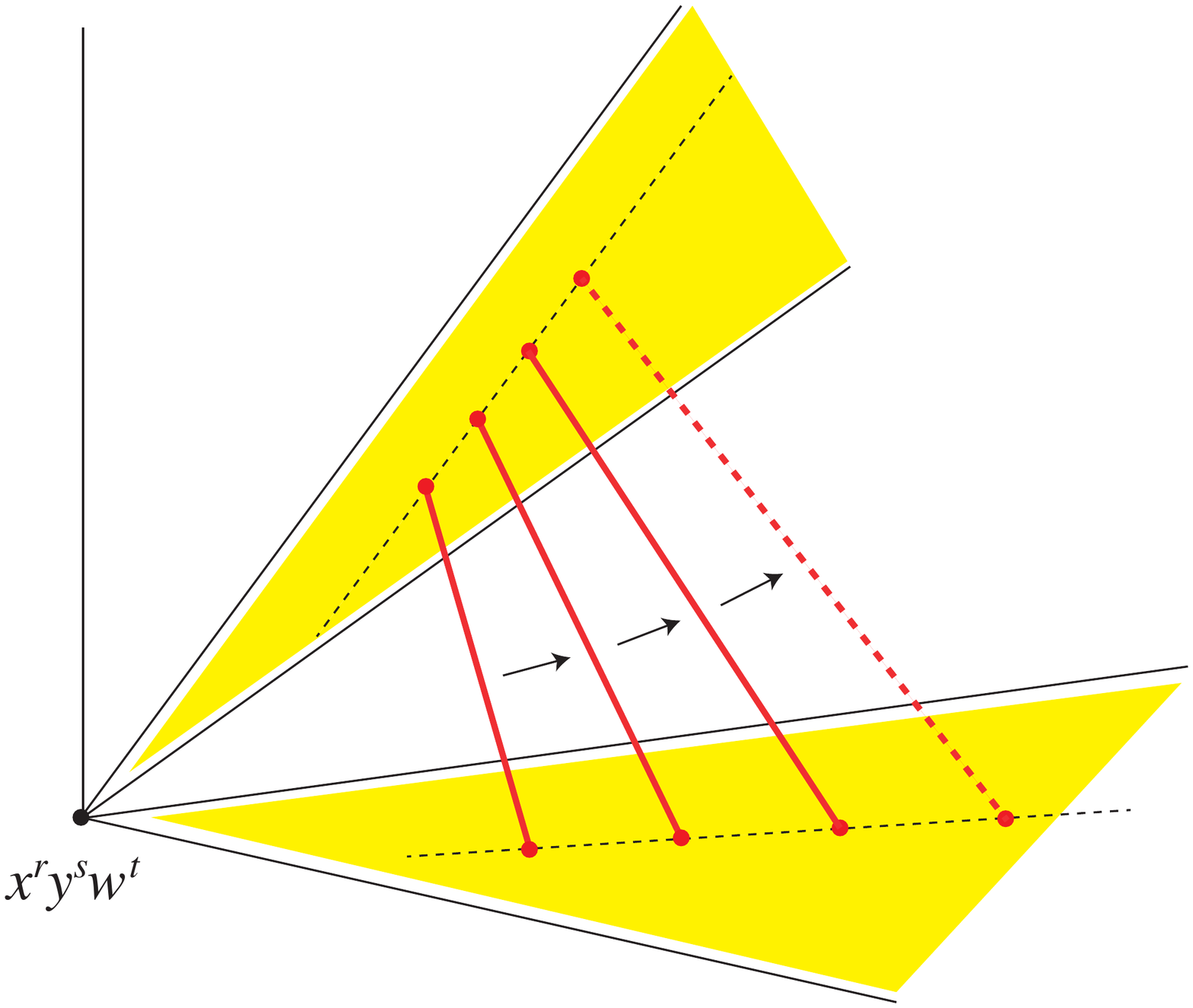} 
\vskip-.4cm

\centerline{{\it Figure} 1: Transformation of the significant edge of the Newton polyhedron.}\medskip\medskip


\st This example follows the same pattern as the first example. Again, the residual order increases only under blowup (1) and the remaining blowups are performed to rebuild the starting equation. Under the last blowup of sequence (2), the residual order decreases by $\frac{p-1}{2}$ and the initial form of $G$ changes to $x^{d'+\frac{p^2+1}{2}}$. The blowup sequences (3) and (4) leave the initial form unchanged, but bring the first term into the form $w^{d'}y^{p^2}\cdot\vt$. Under blowup (5) the residual order decreases by $1$ and $w^{d'}y^{\frac{p^2-1}{2}}$ becomes the new initial form. Again, the last sequence of blowups is only needed to bring the terms of higher degree into the same form as in the original equation. 


Actually, the exponents of both minimal monomials of $F$ increase. The relevant edge of the Newton polyhedra of $F$ and its transforms under blowup connecting the two monomials of the weighted initial form of $F$
are illustrated schematically in Figure 1 by the red segments. The yellow regions represent the $xv$-, respectively $yw$-planes. The segments become longer and move from left to right with each cycle.%


\ignore

\medskip


\section{Analysis of example 2}

\noindent The polynomials of the examples have to satisfy very specific properties so as to produce the increase of the residual order under the first blowup. It is this increase which then gives enough flexibility to obtain by the further blowups a final polynomial having the same shape as at the beginning. We describe below these restrictions for the second example, complementing them by some more general remarks. Similar considerations apply to the first.

\begin{enumerate}

\item Notice that the presentation of the starting polynomial is ``intrinsic''. By this we mean (heuristically speaking) that no coordinate change can simplify the equation or eliminate one of the three indicated monomials. They are all significant. Omitting monomials of the expansion which are multiples of others, the polynomial $f$ has the form
\[z^{p^e}+x^ry^sw^t(w^dy^b+x^av^c),\]
with positive integers $a,b,c,d,r,s,t$ as exponents. Classically, this is (together with the knowledge of the exceptional multiplicities) the entire information one uses and works with when constructing the local invariant for the resolution in zero characteristic. Here, $x,y$ and $w$ are treated as variables defining the exceptional components produced by earlier blowups.

\item The first blowup realizes the increase of the residual order. It passes from $d+\frac{p^2-1}{2}$ to $d+p^2$, hence increases by $\frac{p^2+1}{2}$. For this to happen it is necessary (see \cite{Moh, HP_Characterizing}) that the order of $F$ is a multiple of the order of $f$. The order of $F$ in the second example is $2p^3$. The power $p^3$ as the order of $f$ is convenient so as to permit a sufficiently large increase of the residual order in the first blowup. Recall that by Moh's theorem, the increase is bounded by $p^{e-1}$ when $p^e$ is the order of $f$. In the example, this maximal increase is not realized.

\item The variable $w$ is a {\it silent} variable. By this we understand that the reference point after the blowup variety lies inside the strict transform of the hypersurface $w=0$ (i.e., $w$ does not become invertible in the local ring). Such variables have to appear solely as $p^e$-th powers in order to have an increase of the residual order. The variable $x$ defines the $x$-chart of the blowup (where the reference point is situated), and $y$ and $v$ are {\it relevant} variables, i.e., the transforms of the hypersurfaces $y=0$ and $v=0$ do not contain the reference point (i.e., $y$ and $v$ become invertible in the local ring). In terms of coordinate transformations, this corresponds to the translations $y\mapsto y+1$ and $v\mapsto v+1$. For the increase of the residual order, the second translation $v\mapsto v+1$ is actually irrelevant. The exceptional multiplicities of $x$, $y$ and $v$ are $\frac{p-1}{2}p^2$, $\frac{p^3+1}{2}$, and $0$, respectively. The sum of their residues modulo $p^e$ is $\leq 2p^3$, another necessary condition for the increase of the residual order \cite{HP_Characterizing}.

\item The (weighted) initial form of the starting polynomial $F$ consists of the two monomials $w^dy^{\frac{p^2-1}{2}}$ and $x^{d+\frac{p^2+1}{2}}v^\qq$. The degree of the first increases, whereas the second one decreases.

\item The equimultiple locus of the starting polynomial $F$ equals [insert]. Hence larger center than the chosen point are admissible. One could try to add large powers in pure variables in order to establish at the beginning an isolated singularity which would hence force a point blowup. It is, however, not probable that these powers can be recovered in a similar form after having run through one cycle.

\item Note that $(z,x,w)$ and $(z,v,w)$ are also admissible centers at the beginning, and the respective blowups would not produce the increase of the residual order. The second choice would be taken when following the construction of the center via iterated coefficient ideals \cite{EH}.

\item The monomial $v^\qq\cdot P$ in the starting equation is indicated as the second term of the factor multiplying $w^dy^{\frac{p^2-1}{2}}$. This information is needed to show that no cancellation can occur (which could happen for terms of smaller degree). We are grateful to H.~Kawanoue for pointing out this fact.

\item It seems that taking into account also orders along regular curves as in \cite{HP_Surfaces} for the definition of a suitable resolution invariant does not really help in this example. This is due to the fact that, at each step of the sequence, the initial forms are  monomials.

\item Hironaka chooses for his proof of embedded surface resolution admissible centers of maximal dimension \cite {Hironaka_Bowdoin}. They are not defined by the stratum of maximal values of an invariant but rather ad hoc. His invariant $\delta$ [check which one denotes the vertical height] would increase if centers are not chosen maximal. There is no statement about the upper semicontinuity of the invariant (this is actually not needed, as the centers are chosen directly). And for the induction argument it suffices to know that the invariant takes only finitely many values. [Check whether Abhyankar uses maximal centers, see Cutkosky Skeleton.]

\item {}[?] Our example does not only exhibit an indefinite increase of the residual order. It rather seems that any reasonable invariant attached to $F$ does increase.
[One should try to calculate Cossart's $\omega$ given as the order of the logarithmic Jacobian ideal.] Also, it is not clear how the differential closure of $f$ behaves under the given sequence of blowups, and whether it is able to reveal some additional significant information.

\item When passing from the starting $F$ to the final $F$ after one cycle, the exceptional multiplicities of $x$ and $y$ have remained constant, as well as the exponent of $y$ in the terms of the parenthesis. All other exponents change according to the shift $d\mapsto d'=d+\frac{p^2-p}{2}$.

\end{enumerate}

\section{Conclusions}

\noindent What can we learn from these examples? Do they hint at the existence of a counter-example to the resolution of singularities in positive characteristic, or do they suggest new methods of how to attack and advance in this problem? We are not able to answer these questions to our satisfaction, nevertheless we wish to make a couple of comments which may shed some light on possible answers.

\begin{enumerate}

\item In a (decent) counter-example to resolution in positive characteristic, any choice of an (admissible) center should produce a cycle under blowup (or, at least, a definite deterioration of the singularities). One standard requirement is that the center is contained in the singular locus of the scheme to be resolved (outside the singular locus, there is no need to modify the scheme). Since Zariski, one even requires that the center is regular and contained in the equimultiple locus. This is a reasonable assumption. But even if we start with a singularity whose equimultiple locus is reduced to a single point (e.g., by adding high powers in the variables to the defining polynomial), the equimultiple locus after blowup may be positive dimensional \cite{Ha_Obstacles}. It seems hard to construct a sequence of blowups where at each step the equimultiple locus is zero dimensional and where, at the end, a cycle is produced. 

\includegraphics[width=0.8\textwidth]{edgewithoutzero.pdf}

\centerline{Figure 1.}\bigskip

\item The increase of the residual order in the first blowup is just one aspect of the examples. Actually, the whole initial form of the coefficient ideal of $f$, say, of $F$, has increasing exponents. The transformation of the relevant edges of the Newton polyhedron connecting $x^ry^sw^t\cdot w^dy^b$ with $x^ry^sw^t\cdot x^av^c$ under blowup is depicted in Figure 1.

\item In the resolution algorithm over fields of zero characteristic, positive dimensional centers are only chosen when some coefficient ideal is trivial, say reduced to $0$, or when the so-called monomial case is reached. Choosing in the first case a smaller center than the prescribed one does no harm since the resolution invariant stays constant. But it does not help to improve the singularities. In the second, monomial case, however, the resolution invariant (which is then defined differently from the first case) does increase if the chosen center is too small. This is the phenomenon sketched in the example at the beginning of the note. It produces a cycle with shifted exponent, but it is very easy to circumvent the incident to happen since the correct centers can be prescribed combinatorially. 

\item In our examples, we are not in the monomial case. This suggests to factor exceptional components from $F$ and to look for an improvement in the residual part. Nevertheless, even when doing so, the situation gets worse. This is forced by the necessity of changing the hypersurface of weak maximal contact after the first blowup, say, by cleaning.

\item It would be interesting to test the examples against the modern approaches to resolution in positive characteristic. E.g., to determine the behaviour of the invariant $\omega$ used by Cossart-Piltant and Cossart-Jannsen-Saito, or to compute the differential closure of $F$ as suggested by Kawanoue-Matsuki, respectively the associated elimination algebras of Villamayor and his collaborators [cite ...]. Also the recent proposals of Hironaka should be revisited in the context of the examples [cite ...].

\item The main issue which emerges from the examples is the question of how to detect the correct choice of centers. Of course, in the examples, this can be done ad-hoc, but there seems to be no systematic procedure known of how to do this. This is due, in part, to our limited understanding of characteristic $p$ phenomena under blowup.

\ignore

\item The main conclusion concerns implications and restrictions for prospective proofs of embedded resolution in positive characteristic: Most proofs aim at defining a local upper semicontinuos invariant $\inv_\xi X$ which is associated to each point $\xi$ of the scheme $X$ to be resolved. The invariant stratifies $X$ into locally closed strata. The stratum of points where the invariant attains its maximal value is closed, and the center of blowup is chosen equal to (or at least included in) this stratum. Moreover, the local invariant decreases at each point of the exceptional component created by the blowup in the appropriate center. 

In addition to the above properties of a resolution invariant, it is often chosen in practice  in terms of a (suitably defined) Newton polyhedron, or, said differently, in term of the monomials of the expansion of $f$ whose exponents are minimal with respect to the componentswise order (these choices have to be made in a manner not to depend on any subjective choice of coordinates). Monomials of the expansion of $f$ which are multiples of other monomials are often ignored in the definition of the invariant.

This essential part of starting polynomial in example 2 is of the form \medskip

\hskip 2cm $f=z^{p^e}+x^ry^sw^t(w^dy^b+x^av^c)$\medskip

\noindent and the essential part of the final polynomial is of the form \medskip

\hskip2cm $f'=z^{p^e}+x^ry^sw^{t'}(w^{d'}y^b+x^{a'}v^{c'})$, \medskip

\noindent where $t=kp^3-d$ for some $k$ and where $t'$, $a'$, $c'$ and $d'$ are obtained from $a$, $c$ and $d$ by shifts [not all the same shift, and $t'$ changes differently via changes of $k$]. 

Recall here that in characteristic $0$, the characteristic zero resolution invariant as in \cite{EH} does not  take into account multiplication by monomials (in the first stage of the resolution, until the order of the residual factor $I$ becomes too small). But in the present examples, the polynomials $F^{(i)}$, respectively their initial forms, change more drastically, extending the distance of their exponents, which, themselves, increase.

\item{} [?] On an even more speculative level, the examples show that there are little chances for finding a {\it conveniently hard} polyhedral game whose solution by a winning strategy would imply local embedded resolution (i.e., local uniformization). Again, the problem lies in the required upper semicontinuity of any attached invariant.


\end{enumerate}

\recognize
\providecommand{\bysame}{\leavevmode\hbox to3em{\hrulefill}\thinspace}
\providecommand{\MR}{\relax\ifhmode\unskip\space\fi MR }
\providecommand{\MRhref}[2]{%
  \href{http://www.ams.org/mathscinet-getitem?mr=#1}{#2}
}
\providecommand{\href}[2]{#2}

\smallskip 
\noindent Faculty of Mathematics\par
\noindent University of Vienna, Austria\par
\noindent herwig.hauser@univie.ac.at\par
\noindent stefan.perlega@univie.ac.at

\end{document}